\documentclass[a4paper, 12pt]{amsart}
\usepackage{amssymb,amsmath,amsthm,a4,cancel}
\numberwithin{equation}{section}
\usepackage{enumerate,color}

\newtheorem{thm}{Theorem}[section]
\newtheorem{lem}[thm]{Lemma}

\newtheorem{prop}[thm]{Proposition}
\theoremstyle{definition}

\newtheorem{defn}[thm]{Definition}

\newcommand{\R}{\mathbb{R}}
\newcommand{\C}{\mathbb{C}}
\newcommand{\D}{\mathbb{D}}

\newcommand{\N}{\mathbb{N}}
\newcommand{\T}{\mathbb{T}}

\begin{document}
\title[Similarity to an isometry of composition operators on $\C^+$]{Similarity to an isometry of composition operators on the half-plane}
\author[S.~J.~Elliott]{Sam Elliott}

\address{Department of Pure Mathematics\\
  University of Leeds\\
  Leeds\\
  LS2~9JT\\
  UK}

\email{samuel@maths.leeds.ac.uk}

\subjclass[2000]{47B33}

\begin{abstract}
Necessary and sufficient conditions are already known in the Hardy spaces of both the disc and the half plane for a composition operator to be an isometry, by Nordgren in the disc \cite{Nordgren} and by Chalendar and Partington in the half plane \cite{Chalendar03}. All the same, conditions for such an operator to be similar to an isometry have taken much longer to find. We present some necessary conditions for general weighted composition operators to be similar to an isometry, and use them to produce a complete characterisation of the rational composition on $H^p(\C^+)$ which have this property.
\end{abstract}
\maketitle
\section{Introduction}
Letting $m$ denote normalised Lebesgue measure on the unit circle $\T$, we recall that for $1\le p<\infty$ the Hardy space $H^p(\D)$ is the collection of analytic function on the disc for which
\[
 \|f\|^p = \sup_{0<r<1}\int_\T |f(rt)|^p dm(t) <\infty.
\]
For an analytic map $\phi:\D\to\D$, the composition operator with symbol $\phi$ is the map
\[
 C_\phi f = f \circ \phi,
\]
which can be defined on a number of spaces.
\begin{defn}
We say an operator $T$ is \textit{similar to an isometry} if there is an invertible operator $S$ and an isometry $U$ such that
\[
 T=S^{-1}US.
\]
\end{defn}

In \cite{Bayart}, Fr\'ed\'eric Bayart proved that a composition operator $C_\phi$ on $H^2(\D)$ is similar to an isometry if and only if $\phi$ is inner and has a fixed point in the disc. This result correlates well with Nordgren's result \cite{Nordgren} that such an operator is a genuine isometry if and only if $\phi$ is inner and fixes $0$.

In the case of the half plane, the condition for being an isometry is somewhat less concrete, but nonetheless has echoes of the case of the disc. In this paper, we use the upper half-plane, which will be more convenient. In this setting, yet again, in order for $C_\phi$ to be isometric it must be inner, which is to say it must map the boundary of the half plane to itself almost everywhere, but must also satisfy the condition
\[
 \int_\T \left|\frac{1+\Phi(z)}{1+z}\right|^2 dm(z) = 1,
\]
where $\Phi$ represents the function which is equivalent to $\phi$, but defined on the disc. This condition is a fairly natural one to arise, as we note that $C_\phi$ on $H^2(\C^+)$ is equivalent to the weighted composition operator
\[
 \frac{1+\Phi(z)}{1+z} C_\Phi
\]
on $H^2(\D)$.

In this paper, we will prove a number of results about similarity to an isometry, both of weighted composition operators on the disc, and composition operators on the half plane, ending with a complete characterisation of the rationally induced composition operators on $H^p(\C^+)$.

\section{Preliminaries}
Let $\C^+$ denote the upper half of the complex plane, that is
\[
 \C^+ = \{z\in\C:\Im(z)>0\}.
\]
We recall from, for example \cite{Hoffman}, that the Hardy space $H^p(\C^+)$ is the collection of analytic functions $f$ on $\C^+$ for which the norm
\[
 \|f\|^p = \sup_{y>0}\int_\R |f(x+iy)|^p dx
\]
is finite. In fact, by extending $H^p$ functions non-tangentially to $\R$, it is in fact the case that $H^p$ is a subspace of $L^p(\R)$, and so in particular $H^2$ is a Hilbert space, with inner product
\[
 \left<f,g\right> = \int_\R f(x)\overline{g(x)}dx.
\]
The disc and half plane are conformally equivalent via the M\"obius transformation
\[
 J:\D\to\C^+ \qquad z\mapsto i\left(\frac{1-z}{1+z}\right),
\]
which also gives us the following Lemma. For $\phi:\C^+\to\C^+$, we define a composition operator in the natural way.
\begin{lem}
 If $\phi:\C^+\rightarrow\C^+$ is an analytic self-map of the upper half-plane, then the composition operator $C_\phi:H^p(\C^+)\rightarrow H^p(\C^+)$ (similarly $L^p(\R)\rightarrow L^p(\R)$) is unitarily equivalent to the weighted composition operator $L_\Phi:H^p(\D)\rightarrow H^p(\D)$ (similarly $L^p(\T)\rightarrow L^p(\T)$), given by
\[
 (L_\Phi f)(z) = \left(\frac{1+\Phi(z)}{1+z}\right)^{2/p} C_\Phi f(z) \text{,}
\]
where $\Phi=J^{-1}\circ\phi\circ J$.
\end{lem}

\section{Some necessary conditions}
We begin with the following theorem of B\'ela Sz\H{o}kefalvi-Nagy \cite{Nagy}.
\begin{thm}\label{NagyCondition}
 An operator $T$ on a Hilbert space $\mathcal{H}$ is similar to an isometry if and only if there is some constant $k>0$ such that
\[
 \frac{1}{k}\|x\| \le \|T^nx\| \le k\|x\|,
\]
for all $x\in\mathcal{H}$ and all $n\in\N$.
\end{thm}

The property is certainly also necessary for an operator on a general Banach space to be similar to an isometry, since if $T=S^{-1}US$ then for all $x$ of norm $1$ and all $n\in\N$, we have
\[
 \left\|T^nx\right\| = \left\|S^{-1}U^nSx\right\| \le \underbrace{\|S^{-1}\|\cdot\cancelto{1}{\|U^n\|}\cdot\|S\|}_{=k}\cdot\|x\|.
\]
Similarly, since $S^{-1}$ is bounded, $S$ must be bounded below, and so the lower bound follows in the same fashion.

It will be helpful to consider weighted composition operators on the disc instead of composition operators on the half plane. The following proposition follows a similar proof to that of Bayart in \cite{Bayart}.
\begin{prop}\label{Chap5:Prop1}
 Let $W_{\phi,\psi}$ be a weighted composition operator on the disc such that $\phi$ is rational. Then if $W_{\phi,\psi}$ is similar to an isometry on some $H^p(\D)$, $\phi$ must be inner.
\end{prop}
\begin{proof}
 Let $W_{\phi,\psi}$ be such an operator, namely
\[
 W_{\phi,\psi} f(z) = \psi(z) f \circ \phi(z).
\]
Suppose that $C_\phi$ is similar to an isometry, and hence that there is some $k>0$ such that for all $f\in H^p(\D)$ and all $n\in\N$
\begin{equation}\label{Chap5:Eqn1}
 \frac{1}{k}\|f\| \le \|W_{\phi,\psi}^n f\|.
\end{equation}
Let $B$ be a Borel set on $\T$, and let $B'= \phi^{-1}(B)\cap\T$. Define a function $f\in H^p(\D)$ such that
\[
 |f|=
\begin{cases}
 1 & \text{on $B$} \\
 1/2 & \text{on $\T\setminus B$}
\end{cases},
\]
which must exist since $\log|f|$ is a $p$-integrable function, and so we can construct $f$ by the theory of outer functions. Then $\|f^j\|^p \longrightarrow m(B)$ as $j\rightarrow\infty$, and
\[
 \|W_{\phi,\psi}f^j\|^p = \int_{B'} |\psi(z)|^p|f^j\circ\phi(z)|^p dm(z) + \int_{\T\setminus B'} |\psi(z)|^p\underbrace{|f^j\circ\phi(z)|^p}_{=(1/2)^{p\cdot j}} dm(z).
\]
So, as $j\rightarrow \infty$, the second term on the right hand side of the equation tends to zero, and as such,
\[
 \|W_{\phi,\psi}f^j\|^p \rightarrow \int_{B'} |\psi(z)|^p dm(z).
\]
Now taking limits in (\ref{Chap5:Eqn1}), we get
\begin{equation}\label{Chap5:Eqn2}
 \left(\frac{1}{k^p}\right)m(B)\le \int_{B'} |\psi(z)|^p dm(z).
\end{equation}
By Theorem \ref{NagyCondition} and the comments after it, we have that if $W_{\phi,\psi}$ is similar to an isometry then there is some $k\ge0$ such that
\[
 \frac{1}{k}\|f\|\le\|W_{\phi,\psi}^n f\|\le k\|f\|
\]
for all $n\in\N$. Suppose $\phi$ is not inner. Then since $\phi$ is rational, $|\phi|=1$ in at most finitely many places on $\T$, that is
\[
 m\left(\left\{ z\in\T : |\phi(z)|<1\right\}\right)=1.
\]
Set $B=\{z\in\T:|\phi(z)|<1\}$, and $B'=\phi^{-1}(B)\cap\T$. Now by (\ref{Chap5:Eqn2}) we have
\[
 \frac{1}{k^p}m(B)\le\int_{B'}|\psi(z)|^p dm(z),
\]
and hence, since $m(B)=1$,
\[
 \int_{B'} |\psi(z)|^p dm \ge \frac{1}{k^p}.
\]
In particular, therefore, $m(B')>0$. But $B\cap B'=\emptyset$, since $\phi(B)\subset\D$, so
\[
m(B\cup B')=m(B)+m(B')=1+m(B')>1,
\]
which is a contradiction, since the measure of the entire circle is $1$.
\end{proof}
A number of slight extensions to this theorem exist, which we include for completeness. Firstly, it's easy to see that the proof of Theorem \ref{Chap5:Prop1} we used only the fact that $C_\phi$ is bounded below, and not the much stronger condition of being similar to an isometry, so the Theorem will also hold for this more general class of operators.

The first theorems concerning similar to an isometry of composition operator on the disc were proved by Jaoua in \cite{Jaoua}, and dealt uniquely with maps which were analytic on a neighbourhood of the disc. We aim to prove results in the same direction for weighted composition operators, though we can deal with a slightly more general case. We begin with following Lemma.
\begin{lem}\label{Chap5:Lem1}
 Let $\phi$ be analytic on an open neighbourhood $\Omega$ of $\overline{\D}\setminus\{a\}$ for some $a\in\T$ such that $\phi(\D)\subseteq\D$. Then either
\begin{align*}
m\left(\{z\in\T:|\phi(z)|=1\}\right) & = 0,\text{ or}\\
m\left(\{z\in\T:|\phi(z)|=1\}\right) & = 1.
\end{align*}
\end{lem}
\begin{proof}
 Set
\begin{equation}\label{Chap5:Eqn3}
K_\phi = \{z\in\T : |\phi(z)|=1\}.
\end{equation}
Suppose that $m(K_\phi)\neq0,1$, then certainly both $K_\phi$ and $\T\setminus K_\phi$ are non-empty. Take $a\in\T\setminus K_\phi$ and consider the M\"obius map $F$ which maps $\D\to\C^+$ and $a\mapsto \infty$ by rotation of the Riemann sphere.
We denote by $\widetilde{\phi}$ the map on $\C^+$ which is equivalent to $\phi$ via $F$. Now consider the set
\[
 \mathfrak{A}_\phi = \{z\in\R : \widetilde{\phi}(w)\in\R \text{ on a neighbourhood (in $\R$) of $z$}\}.
\]
This set is clearly open in $\R$. Moreover, consider some sequence of points in $\mathfrak{A}_\phi$, $(z_n)_{n\in\N}\to z_0$. Since $\widetilde{\phi}$ is analytic on some neighbourhood of $\C^+$, $\widetilde{\phi}$ is certainly continuous on $\R$, so if $\widetilde{\phi}(z_n)\in\R$ for each $n\in\N$, we must have that $\widetilde{\phi}(z_0)\in\R$ also.

Let $\widetilde{\phi}(z)=\sum_{j=0}^\infty a_j(z-z_0)^j$ near $z_0$, and denote $\widetilde{\rho}(z)=\sum_{j=0}^\infty \overline{a_j}(z-z_0)^j$. Since $\widetilde{\phi}$ is real at $(z_n)_{n\in\N}$ and at $z_0$, $\widetilde{\phi}=\widetilde{\rho}$ at each of the $z_n$, and at $z_0$. By the identity theorem, therefore, $\widetilde{\phi}=\widetilde{\rho}$ on some neighbourhood of $z_0$, so
\[
 \widetilde{\phi}(z)=\sum_{j=0}^\infty a_j(z-z_0)^j
\]
on some neighbourhood of $z_0$, with each $a_j\in\R$. As such, $\widetilde{\phi}(w)\in\R$ on some neighbourhood of $z_0$, so $z_0\in\mathfrak{A}_\phi$.

So $\mathfrak{A}_\phi$ is a closed and open subset of $\R$, and is non-empty, since it contains a limit point, thus it must be the whole of $\R$. By inverse M\"obius transformation, therefore, we must have $\mathfrak{B}_\phi\supseteq\T\setminus\{a\}$, where
\[
 \mathfrak{B}_\phi = \left\{z\in\T:|\phi(w)|=1 \text{ on a neighbourhood of $z$}\right\}.
\]
But $\mathfrak{B}_\phi\subseteq K_\phi$ (as defined in Equation \ref{Chap5:Eqn3}), so $m(K_\phi)=1$, which is a contradiciton.
\end{proof}
With this Lemma, we can now extend Proposition \ref{Chap5:Prop1} to give the following.
\begin{prop}
 Let $\phi$ be analytic on some neighbourhood $\Omega$ of $\overline{\D}\setminus\{a\}$ for some $a\in\T$ such that $\phi(\D)\subseteq\D$. Then if $W_{\phi,\psi}$ is similar to an isometry on some $H^p(\D)$ for some weight $\psi$, $\phi$ must be inner.
\end{prop}
\begin{proof}
The proof is identical to that of Proposition \ref{Chap5:Prop1}, making use of Lemma \ref{Chap5:Lem1} to give us the fact that $\phi$ must either be inner, or map at most a set of measure $0$ onto the circle.
\end{proof}
Since composition operators on the half plane are equivalent to weighted composition operators on the disc, it's clear that we can make the same claims about such operators, under similar assumptions. It is a special case of Proposition \ref{Chap5:Prop1}, for example, that all rationally induced composition operators on the half plane which are similar to an isometry must have a symbol which is inner. With this in mind, we can now use our results to help characterise similarity to an isometry on the half-plane.

\section[Results for composition operators on $\C^+$]{Results for composition operators on the Half plane}
We will be making use of the following from P\'olya and Szeg\H o's book \cite{PolyaSzegoBook}, page 79.
\begin{lem}\label{PolyaSzego}
 The rational maps of $\R$ for which
\[
 \int_\R f(r(t)) dm(t) = \int_\R f(t) dm(t)
\]
for all $f$ for which the right hand side exists are precisely those maps of the form
\[
 r(z) = \pm\left(z + \alpha + \sum_{i=1}^n \frac{\mu_i}{z-\gamma_i}\right),
\]
where $\alpha\in\R$, the $\gamma_i$ are distinct with $\gamma_i\in\R$ for each $i$, and $\mu_i<0$ for each $i$.
\end{lem}
We note that if $r(\C^+)\subseteq \C^+$, then certainly $r$ being of the above form will imply $C_r$ is an isometry on $L^p(\R)$ and hence $H^p(\C^+)$ by applying Lemma \ref{PolyaSzego} to $|f|^p$. 

We are now in a position to prove our main theorem: a classification of the rationally induced composition operators on $H^p(\C^+)$ which are similar to an isometry.
\begin{thm}
 If $r$ is a rational self map of $\C^+$, the following are equivalent:
\begin{enumerate}
 \item the composition operator $C_r$ is an isometry on some (or equivalently every) $H^p(\C^+)$.
 \item the composition operator $C_r$ is similar to an isometry on some (or equivalently every) $H^p(\C^+)$.
 \item $r$ is inner and $\lim_{|z|\to\infty} r(z)/z=1$.
\end{enumerate}
\end{thm}
\begin{proof}
 That $(1)\Rightarrow(2)$ is of course trivial, and we have already seen as a special case of Proposition \ref{Chap5:Prop1} that $(2)$ implies $r$ is inner. By Corollary 3.5 in \cite{ElliottJury}, we see that the second condition in $(3)$ is equivalent to the statement that $C_r$ has norm $1$. Since $\|C_r^n\|=\|C_r\|^n$ for composition operators on the half plane, again by Corollary 3.5 in \cite{ElliottJury}, if $\|C_r\|\neq1$, we would have a sequence $f_n\in H^p(\C^+)$ with each $f_n$ having norm $1$ such that
\[
 \|C_r^n f_n\| \to \infty
\]
or
\[
 \|C_r^n f_n\| \to 0
\]
as $n\to\infty$, which would contradict the Sz\H{o}kefalvi-Nagy condition (Theorem \ref{NagyCondition}). As such, $(2)\Rightarrow(3)$.

To get $(3)\Rightarrow(1)$ we need to show that all rational functions satisfying $(3)$ are of the form given in Lemma \ref{PolyaSzego}. In order to satify condition $(3)$, $r$ must be of the form
\[
 r(z) = \frac{z^n + a_{n-1} z^{n-1} + \ldots + a_0}{z^{n-1}+b_{n-2}z^{n-2}+\ldots+b_0},
\]
with the $a_i$ and $b_i$ all real, in order that $r$ be inner. As such, $r$ must have a partial fraction expansion of the form
\[
 r(z) = z + \alpha + \sum_{i=1}^n \frac{\mu_i}{(z-\gamma_i)^{m_i}},
\]
for some real $\alpha, \mu_i$ and $\gamma_i$, and natural numbers $m_i$. It remains only to show that none of the roots is repeated (that is, each $m_i$ is $1$), and that each $\mu_i$ is negative. In fact, both facts follow entirely from complex number arguments.

Suppose $r$ has an $n$-fold ($n\ge2$) repeated root $\gamma$, and suppose without loss of generality that this root is $0$. Then near $0$, we will have
\[
 r(z) \approx \frac{\mu}{z^n},
\]
for some $\mu$. If $\mu$ is positive, set
\[
z_k=ke^\frac{i\pi}{2n},
\]
for $k>0$. Then for $k$ sufficiently small,
\[
 r(z_k)\approx \frac{\mu}{(ke^{i\pi/2n})^n}=\left(\frac{\mu}{k^n}\right)e^{-\frac{\pi}{2}i}.
\]
Since $\mu$ is positive and $n\ge2$, each $z_k\in\C^+$, but for $k$ small enough, $r(z_k)$ will have negative imaginary part, but $r$ is a self-map of $\C^+$, which is a contradiction. If $\mu$ is negative, the same trick will work with $z_k=ke^{3i\pi/2n}$.

So all the roots of $r$ must be distinct. To prove that each $\mu_i$ is in fact negative, we again assume without loss of generality that we have a root at $0$. Near $0$, we have
\[
 r(z) \approx \frac{\mu}{z},
\]
so consider the points
\[
 y_k=ke^\frac{i\pi}{2}.
\]
For $k$ sufficiently small,
\[
 r(y_k) \approx \frac{\mu}{ke^\frac{i\pi}{2}} = \left(\frac{\mu}{k}\right)e^{-\frac{\pi}{2}i},
\]
so if $\mu>0$, we have some $r(y_k)$ with a negative imaginary part, but of course each $y_k\in\C^+$, which is again a contradiction to $r$ being a self-map of $\C^+$.
\end{proof}
Although this result holds only for rational functions, it seems not unreasonable to conjecture that, much as in the disc, all composition operators similar to an isometry on $H^2(\C^+)$ ought to have an inner symbol. Indeed, it seems likely that this condition will extend to all $p$ in both the disc and half-plane cases. The equivalence of all three conditions for general maps might be a little too much to ask, but would certainly be an interesting result if it were true, and the equivalence of $(2)$ and $(3)$ is certainly a possibility.


\end{document}